\documentclass{article}

\usepackage[english]{babel}

\usepackage[letterpaper,top=2cm,bottom=2cm,left=3.5cm,right=3.5cm,marginparwidth=1.75cm]{geometry}

\usepackage[utf8]{inputenc}
\usepackage[OT2]{fontenc}
\usepackage[T1]{fontenc}

\usepackage{palatino}

\usepackage{amsmath}
\usepackage{amsfonts}
\usepackage{mathrsfs}
\usepackage{amssymb}
\usepackage{stmaryrd}
\usepackage{chemfig}
\usepackage{mathtools}
\usepackage{amsthm}
\usepackage{wasysym}
\usepackage{xcolor}
\usepackage[colorlinks=true, allcolors=blue]{hyperref}
\usepackage{tikz-cd}
\usetikzlibrary{decorations.markings}
\usetikzlibrary{shapes,positioning,intersections,quotes}

\theoremstyle{plain}
\newtheorem{theo}{Theorem}[section]
\newtheorem{prop}[theo]{Proposition}
\newtheorem{coro}[theo]{Corollary}
\newtheorem{lemma}[theo]{Lemma}

\theoremstyle{definition}

\newtheorem{madef}[theo]{Definition}
\newtheorem{nota}[theo]{Notation}

\newtheorem*{nota*}{Notation}
\newtheorem{rmk}[theo]{Remark}
\newtheorem{ex}[theo]{Example}

\renewcommand{\L}{\mathcal{L}}
\newcommand{\CH}{\underline{\mathrm{CH}}}
\newcommand{\Q}{\mathbb{Q}}
\newcommand{\I}{\mathcal{I}}
\newcommand{\J}{\mathcal{J}}
\newcommand{\M}{\mathrm{M}}
\renewcommand{\H}{\underline{\mathrm{H}}}
\renewcommand{\S}{\mathbb{S}}
\newcommand{\Des}{\mathrm{Des}}
\newcommand{\des}{\mathrm{des}}
\newcommand{\A}{A}
\newcommand{\B}{\mathrm{B}}
\newcommand{\Card}{\mathrm{Card}}
\newcommand{\NM}{\mathsf{NM}}
\newcommand{\Id}{\mathrm{Id}}
\renewcommand{\U}{\mathrm{U}}
\newcommand{\exc}{\mathrm{exc}}
\newcommand{\Exc}{\mathrm{Exc}}
\newcommand{\Asc}{\mathrm{Asc}}
\newcommand{\asc}{\mathrm{asc}}
\newcommand{\D}{\mathcal{D}}
\newcommand{\Hilb}{\mathrm{Hilb}}
\newcommand{\Le}{\mathrm{L}}
\newcommand{\Scol}{\mathcal{S}}
\newcommand{\Gen}{\mathrm{Gen}}
\newcommand{\Is}{\mathcal{I}_{\mathrm{ss}}}
\newcommand{\Isq}{\mathcal{I}_{\mathrm{ss}}^{(2)}}
\renewcommand{\ss}{\mathrm{ss}}
\newcommand{\s}{g}
\newcommand{\fp}{\mathsf{fp}}
\newcommand{\fz}{\mathsf{fz}}
\newcommand{\K}{\mathbb{K}}
\newcommand{\lt}{\mathrm{lt}}
\renewcommand{\d}{d}
\newcommand{\G}{\mathcal{G}}

\usepackage{scalerel}
\newcommand\hugecup{%
    \scaleobj{1.3}{%
        \bigcup\limits_{S \subset E}
    }
}

\title{An algebraic interpretation of Eulerian polynomials, derangement polynomials, and beyond, via Gröbner methods}
\author{Basile Coron}
\date{}

\begin{document}
\maketitle

\begin{abstract}
Motivated by the question of whether Chow polynomials of matroids have only real roots, this article revisits the known relationship between Eulerian polynomials and the Hilbert series of Chow rings of permutohedral varieties. This is done using a quadratic Gröbner basis associated to a new presentation of those rings, which is obtained by iterating the semi-small decomposition of Chow rings of matroids. This Gröbner basis can also be applied to compute certain principal ideals in these rings, and ultimately reestablish the known connection between derangement polynomials and the Hilbert series of Chow rings for corank 1 uniform matroids. More broadly, this approach enables us to express the Hilbert series of Chow rings for any uniform matroid as polynomials related to the ascent statistics on particular sets of inversion sequences.
\end{abstract}
\section{Introduction}
Matroids provide a unifying axiomatization of various notions of independence in mathematics, such as algebraic independence and linear independence (see Welsh \cite{welsh_matroid_1976}). Chow rings of matroids, introduced by Feichtner and Yuzvinsky in \cite{feichtner_chow_2003}, serve as a combinatorial generalization of Chow rings of wonderful compactifications of hyperplane arrangement complements (De Concini-Procesi \cite{de_concini_wonderful_1995}). Those rings have proved to be very deep and interesting objects, laying the foundations for what might be called the algebraic study of matroids. The Hilbert series of these rings (commonly referred to as Chow polynomials of matroids) have garnered considerable attention, demonstrating numerous properties such as palindromicity, unimodality, log-concavity (Adiprasito-Huh-Katz \cite{Adiprasito_2018}), and $\gamma$-positivity (Ferroni-Matherne-Stevens-Vecchi \cite{FERRONI_2024}). For definitions and discussions of these properties and their appearances in enumerative combinatorics, we refer to Bränden \cite{Bona_2015} Chapter 7. Chow rings of matroids are also known to be Koszul (Mastroeni-McCullough \cite{Mastroeni_2021}) and so they satisfy the so-called Koszul inequalities (see Polishchuk-Positselski \cite{Polischuk_2005} for a reference on Koszulness). A central property which is known to imply all the above, provided palindromicity, is that of real-rootedness (see again Bränden \cite{Bona_2015} Chapter 7), which simply means that the polynomial has only real roots. It is a widely open conjecture that all Chow polynomials of matroids are real-rooted. To the best of the author's knowledge, the only matroids for which this conjecture has been confirmed (modulo direct computations) are the simplest matroids (a.k.a the boolean matroids) and the second simplest matroids (a.k.a the corank 1 uniform matroids). The only reason we can establish real-rootedness in these cases is that, fortuitously, for these matroids the Chow polynomials coincide with very well-studied and pervasive polynomials known as the Eulerian polynomials (see \cite{Pet_2015} for a textbook reference), and the derangement polynomials (see Hameister-Rao-Simpson \cite{Hameister_2018} for a proof). The real-rootedness of Eulerian polynomials is a classic result attributed to Fröbenius \cite{Frobenius_1910}, and the real-rootedness of derangement polynomials is due to Gustafsson-Solus \cite{Gustafsson_2020}. Modern techniques allow us to prove both results in a similar manner, which we briefly outline. For any $n\geq 2$, an $n$-inversion sequence $e$ consists of integers $e_1, \ldots, e_n$ such that for all $i\leq n$ one has $0 \leq e_i \leq i-1$. We denote by $\I_n$ the set of such sequences and by $\D_n$ the set of $n$-inversion sequences which do not end at zero and with no consecutive zeros. If we denote by $\asc(e)$ the number of $i$'s such that $e$ increases between $i$ and $i+1$, a possible definition of Eulerian polynomials $A_n(t)$ and derangement polynomials $d_n(t)$ is 
\begin{equation*}
    A_n(t) = \sum_{e\in \I_n} t^{\asc(e)}\,\,\,, \,\,\, d_n(t) = \sum_{e\in \D_n} t^{\asc(e)}.
\end{equation*}
One can partition $\I_n$ according to the last value of the inversion sequence, which gives two sequences of polynomials 
\begin{equation*}
    A^i_n(t) = \sum_{\substack{e\in \I_n\\ e_n = i}} t^{\asc(e)} \,\,\,, \,\,\, d^i_n(t) = \sum_{\substack{e\in \D_n\\ e_n = i}} t^{\asc(e)}.
\end{equation*}
We can then demonstrate that these two sequences of polynomials are both interlacing (meaning the roots of these polynomials are intertwined in a specific manner) by establishing a recurrence relation for them (see Bränden \cite{Bona_2015} Chapter 7 for more details). Those two interlacing results imply the real-rootedness of $A_n$ and $d_n$ respectively, since a sum of interlacing polynomials is real-rooted. 

A first goal of this article is to reprove the equality between Eulerian polynomials and Chow polynomials of boolean matroids in a way which we think could open the door to a generalization of the above strategy, by using the theory of Gröbner bases. In summary, we will establish a bijection between $n$-inversion sequences and the normal monomials of a quadratic Gröbner basis for Chow rings of boolean matroids obtained in a previous article (see \cite{coron_2023}). We then use a similar but different quadratic Gröbner basis in order to compute the Hilbert series of certain ideals in those rings, which reproves the equality between Chow polynomials of uniform matroids in corank $1$ and derangement polynomials. This new quadratic Gröbner basis is analogous to the first one but is associated to a new presentation obtained by iterating the semi-small decomposition for Chow rings of matroids (Braden-Huh-Matherne-Proudfoot-Wang \cite{Braden_2022}). These methods can then be extended to express Chow polynomials of uniform matroids in any corank as polynomials associated with the ascent statistics over certain subsets of inversion sequences, as in the corank 1 case. 
\\  \\
\textbf{Acknowledgement.} The author would like to express his gratitude toward Alex Fink for his generous availability and many interesting conversations. This work was supported by the Engineering and Physical Sciences Research Council [grant number EP/X001229/1].
\section{Chow rings of boolean matroids and Eulerian polynomials}
We start by recalling a few definitions. 
\begin{madef}[Matroid]
A \textit{matroid} $\M$ on a finite set $E$ is a nonempty collection of subsets of $E$, called \textit{flats} of $\M$, which satisfies the following properties:
\begin{enumerate}
    \item The intersection of any two flats is a flat.
    \item For any flat $F$, any element in $E\setminus F$ is contained in exactly one flat that is minimal among the flats strictly containing F.
\end{enumerate}
The set of flats of $\M$ ordered by inclusion will be denoted by $\L(\M)$. A matroid will be called \textit{loopless} if the empty subset of $E$ is a flat. 
\end{madef}
We refer to \cite{welsh_matroid_1976} for a textbook reference on this subject. For any finite set $E$, the matroid having as flats every subset of $E$ will be called the \textit{boolean matroid} on $E$, and will be denoted $\B_{E}$. 
\begin{madef}[Chow ring of a matroid]
Let $\M$ be a loopless matroid on a finite set $E$. The \textit{Chow ring of} $\M$, denoted by $\CH(\M)$, is the quotient algebra
\begin{equation*}
    \Q[x_F, \, F \in \L(\M) \setminus\{\emptyset\}]/(\I + \J)
\end{equation*}
where $\I$ is the ideal generated by the elements 
\begin{equation*}
    \sum_{\substack{F \in \L(\M) \\ F \ni i}}x_F, \, \textrm{ for all } i\in E,
\end{equation*}
and $\J$ is the ideal generated by the elements
\begin{equation*}
    x_{F_1}x_{F_2}, \, \textrm{ for every pair of incomparable flats } F_1, F_2 \in \L(\M) \setminus\{\emptyset\}. 
\end{equation*}
\end{madef}
Chow rings of matroids are graded, by setting the degree of the generators to one. The Hilbert series $\sum_{k}\dim_{\Q} \CH^k(\M) t^k$ will be called the \textit{Chow polynomial of }$\M$ and will be denoted by  $\H_{\M}(t)$. Chow rings of matroids admit another presentation, called the simplicial presentation. 
\begin{prop}[\cite{Pagaria_2023} Theorem 2.9]
For all loopless matroid $\M$ on a finite set $E$, the Chow ring $\CH(\M)$ is isomorphic to the quotient algebra
\begin{equation*}
    \Q[h_F, F \in \L(\M)\setminus\{\emptyset\}]/(\I + \J)
\end{equation*}
where $\I$ is the ideal generated by the elements $h_{\overline{\imath}} \, (i \in E)$ with $\overline{\imath}$ the intersection of the flats containing $i$, and $\J$ is the ideal generated by the elements $(h_F - h_{F_1})(h_F - h_{F_2})$ for all triples $F, F_1, F_2 \in \L(\M) \setminus \{\emptyset\}$ such that $F$ is the intersection of the flats containing both $F_1$ and $F_2$.
\end{prop}
This presentation was first introduced in \cite{Yuzvinsky_2002} and further studied in \cite{Backman_Spencer_Eur_2020}. For any integer $n\geq 1$ we denote $[n]\coloneqq \{1, \ldots, n\}$. 
\begin{madef}[Eulerian polynomial]
For any finite linearly ordered set $E$ of cardinality $n$, let $\S_{E}$ denote the set of bijections from $[n]$ to $E$. For all $\sigma\in \S_E$, an element $i \in [n-1]$ is a \textit{descent} of $\sigma$ if we have the inequality $\sigma(i) > \sigma(i+1)$. We denote by $\Des(\sigma)$ the set of descents of $\sigma$ and by $\des(\sigma)$ the cardinality of $\Des(\sigma)$. The \textit{$n$-th Eulerian polynomial} $\A_n(t)$ is the polynomial associated to the descent statistics on $\S_{[n]}$, that is,
\begin{equation*}
    \A_n(t) \coloneqq \sum_{\sigma \in \S_{[n]}}t^{\des(\sigma)} . 
\end{equation*}
We refer to \cite{Pet_2015} for more details on this topic. Using the fact that the Chow ring of $\B_{[n]}$ is isomorphic to the Chow ring of the toric variety associated to the $n$-th permutohedron (\cite{feichtner_chow_2003} Theorem 3), and the general fact that the Hilbert series of the Chow ring of a toric variety associated to a simple polytope is its $h$-vector (\cite{Danilov_1978} Theorem 10.8), one gets the following result.
\end{madef}
\begin{prop} [\cite{Braden_2022} Remark 2.11] \label{propeulerchow}
For all $n\geq 1$ we have $$\A_n(t) = \H_{\B_{[n]}}(t).$$
\end{prop}
The main objective of this section is to reprove this statement in a more direct way, using the theory of Gröbner bases. Since that theory will be our main tool in this article, we include below some recollections on that subject. We refer to \cite{BW_1993} for more details. \\

Let $\K$ be any field and $V$ some finite set. Given a total order $\vartriangleleft$ on $V$, the associated \textit{degree lexicographic order} $\vartriangleleft_{deglex}$ is the order on the monomials of the free commutative algebra $\K[V]$ defined by setting $ m = x_{v_1}\ldots x_{v_n} \vartriangleleft_{deglex} m' = x_{v'_1}\ldots x_{v'_m}$ (with $v_1 \trianglelefteq \ldots \trianglelefteq v_n \in V $ and $v'_1 \trianglelefteq \ldots \trianglelefteq v'_m \in V$) if $n < m$, or if $n = m$ and $v_1 \vartriangleleft v'_1$,  or if $n=m$, $v_1 = v'_1$ and $m/x_{v_1} \vartriangleleft_{deglex} m'/x_{v'_1}$. Having chosen $\vartriangleleft$, for any nonzero element $f \in \K[V]$ we call \textit{leading term} of $f$ and denote by $\lt(f)$ the greatest monomial for the order $\vartriangleleft_{deglex}$ whose coefficient in the expansion of $f$ is nonzero. For any ideal $\I \subset \K[V]$, a subset $\G \subset \I$ is called a \textit{Gröbner basis} of $\I$ (with respect to $\vartriangleleft$) if every leading term of some element of $\I$ is divisible by the leading term of some element of $\G$. A Gröbner basis is called \textit{quadratic} if it only contains homogeneous elements of degree $2$. A monomial is said to be \textit{normal} with respect to some subset $\G \subset \I$ if it is not divisible by the leading term of some element in $\G$. We have the key proposition. 
\begin{prop}[\cite{BW_1993} Theorem 5.35]\label{propmaingrob}
For any $\G \subset \I$, the set of (cosets of) normal monomials with respect to $\G$ forms a linear basis of $\K[V]/\I$ if and only if $\G$ is a Gröbner basis of $\I$.  
\end{prop}
If $\G \subset \I$ is a Gröbner basis, one has a concrete way of expressing any element $f \in \K[V]$ as a linear combination of normal monomials modulo $\I$, via the following procedure. If the monomials with nonzero coefficients in $f$ are all normal then there is nothing to do. Otherwise, consider $m$  the greatest of those monomials which is not normal. By definition of normality, that monomial is divisible by the leading term of some element $R = \lt(R) - \sum \textrm{lower terms}\in \G.$ Inside $m$ one can replace $\lt(R)$ by $\sum \textrm{lower terms}$ which replaces $f$ by another element $f' \in \K[V]$. Repeat this process until there are no nonnormal monomials left. The end result of that procedure (which terminates since $\vartriangleleft_{deglex}$ is a well-order) will be called the \textit{rewriting} of $f$. Gröbner bases also provide a way to compute ideals in the corresponding quotient ring. For instance for any $f\in \K[V]$, the principal ideal $(f) \subset \K[V]/\I$ is spanned by the rewriting of elements of the form $fm$ where $m$ is a normal monomial with respect to $\G$.\\

In view of Proposition \ref{propmaingrob}, to prove Proposition \ref{propeulerchow} it is enough to find a bijection between $\S_{[n]}$ and the normal monomials associated to some Gröbner basis of $\CH(\B_{[n]})$, which sends the number of descents to the degree of the monomial. In \cite{feichtner_chow_2003} the authors computed a Gröbner basis of $\CH(\M)$ (with the $x$-presentation) for any loopless matroid $\M$, but this Gröbner basis is not well suited for our purpose. In \cite{coron_2023} the author described a different Gröbner basis of $\CH(\M)$, this time for the simplicial presentation. This Gröbner basis is valid only for matroids whose lattice of flats are supersolvable (see \cite{Stanley_1972}), which is true for the boolean matroids. The description of the normal monomials associated to that Gröbner basis revolves around the following definition.
\begin{madef}[Initial segment]
Let $I \subset J$ be two nonempty subsets of a linearly ordered finite set $E$. The subset $I$ is an \textit{initial segment} of $J$ if we have the equality
\begin{equation*}
    I = J \cap \{i \, | \, i \in E, i \leq \max I\}. 
\end{equation*}
\end{madef}
In plain English this means that every element in $J$ which is below $\max I$ must also be in $I$. 
\begin{nota}
For every linearly ordered set $E$ we denote by $\Card(E)$ the cardinality of $E$. If $S$ is a subset of $E$ and $e$ is some element of $E$ we denote 
\begin{equation*}
    S_{\leq e} \coloneqq \{s \in S \, | \, s \leq e\},
\end{equation*}
and similarly for $\geq, <$ and $>$.
\end{nota}
The following proposition summarizes the results of \cite{coron_2023} applied to the boolean matroids. 
\begin{prop}[\cite{coron_2023} Proposition 3.2, Theorem 3.6]\label{propgbasish}
For any finite linearly ordered set $E$, the Chow ring $\CH(\B_{E})$ with its simplicial presentation admits a quadratic Gröbner basis, whose quadratic normal monomials are 
\begin{equation} \label{eqnormal}
    \hugecup\left\{ h_{S_1}h_{S_2} \, \,  \begin{array}{|l} 
        \Card(S_1) \geq 2 \textrm{ and } \Card(S_2) \geq 2, \\
         S_1 \cup S_2 = S \textrm{ and } S_1 \textrm{ is the minimal initial segment of } S \textrm{ satisfying that property,} \\
        \Card(S\setminus S_2) \geq 2 \textrm{ unless } S_2 = S\setminus \{\max S\}
    \end{array} \right\}
\end{equation}
\end{prop}
Since the Gröbner basis is quadratic, a monomial is normal with respect to that Gröbner basis if and only if every quadratic submonomial of that monomial is normal. We will denote by $\NM(E)$ the set of those normal monomials. An important feature of $\NM(E)$ is that it is square free. This means that a normal monomial can be described by the set of the subsets indexing its generators.
\begin{ex}
For $E = [3]$ we get the normal monomials $1, h_{12}, h_{13}, h_{23}, h_{123}, h_{12}h_{123} \in \CH(\B_{[3]})$.
\end{ex}
For all $\sigma\in \S_E$ and $i \in \Des(\sigma)$ we denote $S(\sigma, i) \coloneqq \{\sigma(j) \, | \,\sigma(j) \leq \sigma(i), j \geq i \}$.
\begin{prop}\label{propbij}
For all finite linearly ordered set $E$, the map $\Psi$ which sends a bijection $\sigma \in \S_{E}$ to the monomial $\prod_{i\in \Des(\sigma)}h_{S(\sigma, i)}$ establishes a bijection between $\S_{E}$ and $\NM(E)$.
\end{prop}
By definition, the grading of $\Psi(\sigma)$ is the number of descents of $\sigma$, and so Proposition \ref{propbij} gives a new proof of Proposition \ref{propeulerchow}. Here are some examples of images of $\Psi$ below. 
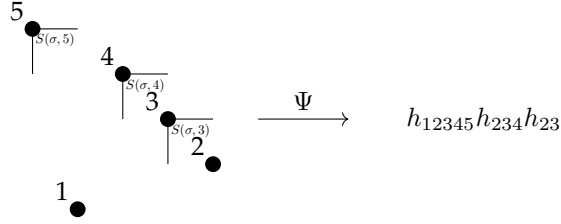
\begin{ex}
\begin{equation*}
\begin{tikzpicture}[scale=0.60, baseline=14px]
    \fill[black] (0,4) circle (5pt)
               (1,0) circle (5pt)
               (2,3) circle (5pt)
               (3,2) circle (5pt)
               (4,1) circle (5pt);
    \node[above left=0pt of {(0,4)}, outer sep=0pt] {5};
    \node[above left=0pt of {(1,0)}, outer sep=0pt] {1};
    \node[above left=0pt of {(2,3)}, outer sep=0pt] {4};
    \node[above left=0pt of {(3,2)}, outer sep=0pt] {3};
    \node[above left=0pt of {(4,1)}, outer sep=0pt] {2};
    \draw[] (0,4) -- (1, 4) node[midway,below, scale = 0.5] {$S(\sigma,5)$} ;
    \draw[] (0,4) -- (0, 3);
    \draw[] (2,3) -- (3, 3) node[midway,below, scale = 0.5] {$S(\sigma,4)$};
    \draw[] (2,3) -- (2, 2);
    \draw[] (3,2) -- (4, 2) node[midway,below, scale = 0.5] {$S(\sigma,3)$};;
    \draw[] (3,2) -- (3, 1);
    \draw[->] (5,2) -- (7,2) node[midway,above] {$\Psi$};
    \node[] (A) at (10,2) {$h_{12345}h_{234}h_{23}$};
\end{tikzpicture}
\end{equation*}
\newline
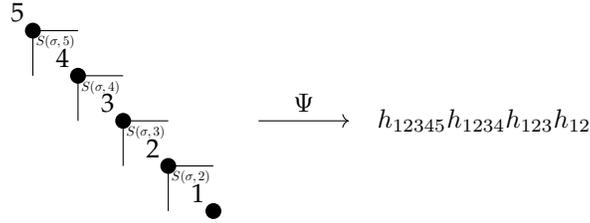
\begin{equation*}
\begin{tikzpicture}[scale=0.60, baseline=14px]
    \fill[black] (0,4) circle (5pt)
               (1,3) circle (5pt)
               (2,2) circle (5pt)
               (3,1) circle (5pt)
               (4,0) circle (5pt);
    \node[above left=0pt of {(0,4)}, outer sep=0pt] {5};
    \node[above left=0pt of {(1,3)}, outer sep=0pt] {4};
    \node[above left=0pt of {(2,2)}, outer sep=0pt] {3};
    \node[above left=0pt of {(3,1)}, outer sep=0pt] {2};
    \node[above left=0pt of {(4,0)}, outer sep=0pt] {1};
    \draw[] (0,4) -- (1, 4) node[midway,below, scale = 0.5] {$S(\sigma,5)$};
    \draw[] (0,4) -- (0, 3);
    \draw[] (1,3) -- (2, 3) node[midway,below, scale = 0.5] {$S(\sigma,4)$};
    \draw[] (1,3) -- (1, 2);
    \draw[] (2,2) -- (3, 2) node[midway,below, scale = 0.5] {$S(\sigma,3)$};
    \draw[] (2,2) -- (2, 1);
    \draw[] (3,1) -- (4, 1) node[midway,below, scale = 0.5] {$S(\sigma,2)$};
    \draw[] (3,1) -- (3, 0);
    \draw[->] (5,2) -- (7,2) node[midway,above] {$\Psi$};
    \node[] (A) at (10,2) {$h_{12345}h_{1234}h_{123}h_{12}$};
\end{tikzpicture}
\end{equation*}
\newline
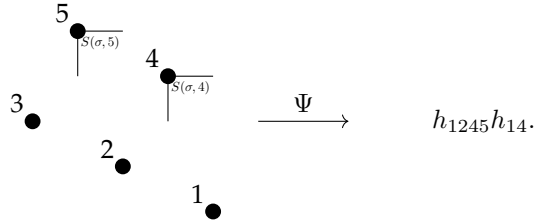
\begin{equation*}
\begin{tikzpicture}[scale=0.60, baseline=14px]
    \fill[black] (0,2) circle (5pt)
               (1,4) circle (5pt)
               (2,1) circle (5pt)
               (3,3) circle (5pt)
               (4,0) circle (5pt);
    \node[above left=0pt of {(0,2)}, outer sep=0pt] {3};
    \node[above left=0pt of {(1,4)}, outer sep=0pt] {5};
    \node[above left=0pt of {(2,1)}, outer sep=0pt] {2};
    \node[above left=0pt of {(3,3)}, outer sep=0pt] {4};
    \node[above left=0pt of {(4,0)}, outer sep=0pt] {1};
    \draw[] (1,4) -- (2, 4) node[midway,below, scale = 0.5] {$S(\sigma,5)$};
    \draw[] (1,4) -- (1, 3);
    \draw[] (3,3) -- (4, 3) node[midway,below, scale = 0.5] {$S(\sigma,4)$};
    \draw[] (3,3) -- (3, 2);
    \draw[->] (5,2) -- (7,2) node[midway,above] {$\Psi$};
    \node[] (A) at (10,2) {$h_{1245}h_{14}.$};
\end{tikzpicture}
\end{equation*}
\end{ex}
\begin{proof}
    Our first step is to prove that for all $\sigma\in \S_{E}$, the monomial $\Psi(\sigma)$ is a normal monomial. Since the Gröbner basis we are using is quadratic, we only need to prove that for any pair $i<j$ of descents of $\sigma$ the quadratic monomial $h_{S(\sigma, i)}h_{S(\sigma,j)}$ is normal. The cardinality requirements in \eqref{eqnormal} are inherent in the definition of $\Psi$. If $p\in E$ is such that $p$ belongs to $S(\sigma, i) \cup S(\sigma, j)$ and $p$ is less or equal than $\max S(\sigma,i) = \sigma(i)$ then $p$ is of the form $\sigma(j')$ for some $j'\geq i$ and so we must have $p \in S(\sigma, i)$. This proves that $S(\sigma, i)$ is an initial segment of $S(\sigma,i) \cup S(\sigma,j)$. Furthermore, $\sigma(i)$ does not belong to $S(\sigma, j)$ and so $S(\sigma, i)\setminus \{\sigma(i)\}\cup S(\sigma,j)$ is strictly contained in $S(\sigma,i) \cup S(\sigma,j)$, which proves the minimality requirement in \eqref{eqnormal}. \\

    We shall now describe the inverse of $\Psi$. Let $m = \prod_{S \in V}h_S$ be a normal monomial in $\CH(\B_{E})$, with $V$ some set of subsets of $E$. We define a binary relation $\vartriangleleft$ on $V$ by setting $S_1 \vartriangleleft S_2$ if $S_1$ is an initial segment of $S_1\cup S_2$ such that $(S_1\setminus\{\max S_1\})\cup S_2$ is strictly contained in $S_1 \cup S_2$. 
    \begin{lemma}
        The binary relation $\vartriangleleft$ defines a strict total order on $V$. 
    \end{lemma}
    \begin{proof}
    If $S_1 \vartriangleleft S_2$ and $S_2 \vartriangleleft S_1$ then both $S_1$ and $S_2$ are initial segments of $S_1 \cup S_2$, which implies that one of those two subsets must be included in the other. In that case removing the maximal element of the smaller subset would not change the union, which is a contradiction. By description \eqref{eqnormal} of the quadratic normal monomials, if $S_1 \neq S_2$ then we must have $S_1 \vartriangleleft S_2$ or $S_2\vartriangleleft S_1$. Let $S_1, S_2, S_3$ be three subsets in $V$ such that we have $S_1\vartriangleleft S_2$ and $S_2 \vartriangleleft S_3$. We are looking to prove $S_1 \vartriangleleft S_3$. We split the proof in  different cases. \\
    
    1) $S_1 \supsetneq S_2 \supsetneq S_3$.\\
    
    In that case $S_1$ is an initial segment of $S_1 \cup S_3 = S_1$. Furthermore, since $(S_1\setminus\{\max S_1\})\cup S_2$ is strictly contained in $S_1 \cup S_2$, the subset $S_2$ does not contain $\max S_1$ and therefore $S_3$ does not either, which implies that $S_1\setminus \{\max S_1\} \cup S_3$ is strictly contained in $S_1 \cup S_3$. In other words, we have $S_1 \vartriangleleft S_3$. \\
    
    2) $S_1 \supsetneq S_2$ and $S_2$ is not comparable with $S_3$ (for $\subset$).\\
    
    In that case we have $\max S_3 > \max S_2$ and $\max S_2$ does not belong to $S_3$. This means that $\max S_2$ is an element smaller than $\max S_3$ which belongs to $S_1 \cup S_3$ and not to $S_3$. This implies that $S_3$ is not an initial segment of $S_1 \cup S_3$ and we must have $S_1 \vartriangleleft S_3.$\\
    
    3) $S_1$ is not comparable with $S_2$ (for $\subset$), and $S_2 \supsetneq S_3$.\\
    
    In that case $\max S_1$ does not belong to $S_2$, and hence does not belong to $S_3$ either. If $\max S_1 < \max S_3$ then as in the case above $S_3$ cannot be an initial segment of $S_1 \cup S_3$, and so $S_1 \vartriangleleft S_3$. On the other hand if $\max S_1 > \max S_3$ then by the fact that $S_3$ is contained in $S_2$ and the fact that $S_1$ is an initial segment of $S_1 \cup S_2$ we get $S_3 \subset S_1$ which also implies $S_1 \vartriangleleft S_3$.\\
    
    4) $S_1$ is not comparable with $S_2$ (for $\subset$), and $S_2$ is not comparable with $S_3$ (for $\subset$).\\
    
    In that case we have $\max S_3 > \max S_2 > \max S_1$ and $\max S_1$ does not belong to $S_2$. Since $S_2$ is an initial segment of $S_2 \cup S_3$ this implies that $\max S_1$ does not belong to $S_3$ either. As in the previous cases this implies that $S_3$ cannot be an initial segment of $S_1\cup S_3$, which implies $S_1 \vartriangleleft S_3$. 
    \end{proof}
    Denote $V = \{S_1 \vartriangleleft \ldots \vartriangleleft S_k\}$. For all $i \leq k$ we define $X_i \coloneqq E_{\leq \max S_i} \setminus (S_i\setminus\{\max S_i\})$. Let $\Phi$ be the map which sends a normal monomial $m$ to the bijection $\Phi(m)$ whose list of elements $(\Phi(m)(1))(\Phi(m)(2))\ldots$ is given by 
    \begin{equation}\label{eqdefphi}
    X_1 (X_2\setminus X_1)(X_3 \setminus (X_1\cup X_2))\ldots(X_k \setminus (X_1\cup \ldots \cup X_{k-1}))([n]\setminus (X_1 \cup \ldots \cup X_k))
    \end{equation}
    where the concatenation of sets means the concatenation of the lists of their elements in increasing order. 

    \begin{lemma}
        The map $\Phi_E$ is the inverse of $\Psi_E$.
    \end{lemma}
    \begin{proof}
    We will prove this statement by induction on the cardinality of $E$. We start with the equality $\Psi_E\circ \Phi_E = \Id_{\NM(E)}$. Let $m$ be a normal monomial with underlying set of subsets $V = \{S_1 \vartriangleleft \ldots \vartriangleleft S_k\}$. Let $e$ be some element in $X_1 \setminus \{\max S_1\}$. The element $e$ is less or equal than $\max S_1$ and does not belong to $S_1$. Since $S_1$ is an initial segment of $S_1\cup S_i$ for all $i>1$ this implies that $e$ does not belong to $S_i$ for all $i>1$. Furthermore, the element $\max S_1$ does not belong either to any $S_i \,(i>1)$ by the assumption that $S_1\setminus\{\max S_1\} \cup S_i$ is strictly contained in $S_1\cup S_i$ for all $i>1$. Consequently, for all $i>1$ we have $X_1 \cap S_i = \emptyset$, and thus one can view the monomial $m' \coloneqq m/h_{S_1}$ as a normal monomial in $\CH(\B_{E\setminus X_1})$. By construction we have $\Phi_E(m) = X_1\Phi_{E\setminus X_1}(m')$. By construction and induction we have $$\Psi_{E}(\Phi_E(m)) = f_{S(\Phi_E(m), i)}\Psi_{E\setminus X_1}(\Phi_{E \setminus X_1}(m')) = f_{S(\Phi_E(m), i)}m'$$ where $i = \Card(X_1)$, and $f_{S(\Phi_E(m), i)}$ is equal to $h_{S(\Phi_E(m), i)}$ if $i$ is a descent of $\Phi_E(m)$ and $1$ otherwise. What is left is to prove that $S(\Phi_E(m), i)$ is equal to $S_1$ and that $i$ is a descent of $\Phi_E(m)$. Note that $\Phi_E(m)(i) = \max S_1$, and so
    \begin{align*}
        S(\Phi_E(m), i) &= \{\Phi_E(m)(j) \, |\, j \geq i, \Phi_E(m)(j) \leq \max S_1\} \\
        &= [\max S_1] \setminus (X_1 \setminus \{\max S_1\}) \\
        &= [\max S_1] \setminus ([\max S_1]\setminus S_1) \\
        &= S_1.
    \end{align*}
    Finally, to prove that $i$ is a descent of $\Phi_{E}(m)$ we must prove that $\min (X_2\setminus X_1)$ is less than $\max S_1$. If $\max S_2 < \max S_1$ then since $S_1$ is an initial segment of $S_1\cup S_2$ we must have $\max S_2 \in S_1$ and so $\max S_2$ belongs to $X_2 \setminus X_1$. This implies $\min (X_2\setminus X_1) \leq \max S_2 < \max S_1$. Otherwise, by description \eqref{eqnormal} of the quadratic normal monomials we must have $\Card(S_1\cup S_2 \setminus S_2) \geq 2$. This means that there exists some element $s$ of $S_1 \setminus S_2$ which is different from $\max S_1$. The element $s$ is less than $\max S_1$ and so is less than $\max S_2$, and since it does not belong to $S_2$ it must belong to $X_2 \setminus X_1.$ This implies $\min (X_2 \setminus X_1) \leq s < \max S_1$. \\

    Let us now prove the equality $\Phi_E\circ \Psi_E = \Id_{\S_{E}}.$ Let $\sigma$ be a bijection from $[n]$ to $E$. By the proof of the fact that $\Psi_{E}(m)$ is a normal monomial which we gave at the very beginning of this proof, one can see that the order on the subsets $S(\sigma, i)\, i \in \Des(\sigma)$ given by the order on $[n]$ is exactly the order $\vartriangleleft$ on the subsets indexing the generators of $\Psi(\sigma)$ that we have defined previously. Let $i_1$ be the smallest descent of $\sigma$. The values $\sigma(i) \, (i \leq i_1)$ are in increasing order and this set of values is exactly $E_{\leq \max S(\sigma, i_1)}\setminus (S(\sigma, i_1)\setminus \{\max S(\sigma, i_1)\})$. In other words we have proved that the first values until $i_0$ in the list describing $\Phi_E\circ \Psi_E(\sigma)$ are the same as in the list describing $\sigma$. Then, the same induction on the cardinality of $E$ as above concludes the proof.   
    \end{proof}

\end{proof}

\section{Chow rings of uniform matroids and derangement polynomials}
\begin{madef}[Uniform matroid]
For any finite set $E$ and any integer $1 \leq k \leq \Card(E)$, the \textit{rank} $k$ \textit{uniform matroid on $E$}, denoted $\U_{k, E}$, is the matroid whose flats are $E$ and all subsets of $E$ of cardinality less or equal than $k-1$. 
\end{madef}
We have a morphism of rings $\CH(\B_{E}) \rightarrow \CH(\U_{k, E})$ defined by sending $h_{S}$ to $h_{S}$ if $S\subset E$ has cardinality less or equal than $k-1$ and sending $h_{S}$ to $h_{E}$ otherwise. One has the following result, which will be central for us. 
\begin{prop}[\cite{Pagaria_2023} Lemma 4.11] \label{propidealchowpoly}
For any finite set $E$ of cardinality $n$, and any integer $1 \leq k \leq n$ we have an isomorphism of $\CH(\B_{E})$-modules 
\begin{equation*}
    (h_{E}^{n-k}) \simeq \CH(\U_{k, E})[-(n-k)]
\end{equation*}
(where $[-]$ denotes the shift in degree). In particular, we have 
\begin{equation*}
    \Hilb((h_{E}^{n-k})) = t^{n-k}\H(\U_{k, E}).
\end{equation*}
\end{prop}
In corank $1$ it turns out that the Hilbert series of $\CH(\U_{n-1, n})$ are (up to a shift) well-known polynomials called ``derangement polynomials'' which we now describe. 
\begin{madef}[Derangement polynomial]
Let $\sigma$ be a permutation of $[n]$. An element $i\in [n]$ is called an \textit{excedance} of $\sigma$ if one has $\sigma(i) > i$. We denote by $\Exc(\sigma)$ the set of excedances of $\sigma$ and by $\exc(\sigma)$ the cardinality of $\Exc(\sigma)$. The permutation $\sigma$ is called a \textit{derangement} if one has $\sigma(i) \neq i$ for all $i \in [n]$. The $n$-th \textit{derangement polynomial}, denoted $d_n(t)$, is the polynomial associated to the excedance statistics on the set of derangements of $\S_{[n]}$, that is, 
\begin{equation*}
    d_n(t) := \sum_{\substack{\sigma \in \S_[n] \\ \sigma \textrm{ derangement}}} t^{\exc(\sigma)}.
\end{equation*}
\end{madef}
To put Eulerian polynomials and derangement polynomials on the same ground we need the following definition. 
\begin{madef}[Inversion sequence]
For any integer $n \geq 1$, an $n$-\textit{inversion sequence} is a sequence of integers $(e_i)_{1\leq i \leq n}$ such that for all $i \leq n$ one has $0 \leq e_i \leq i-1$. We denote by $\I_n$ the set of $n$-inversion sequences. For any $n$-inversion sequence $e$, an element $i \leq n-1$ is called an \textit{ascent} of $e$ if $e_i < e_{i+1}$. We denote by $\Asc(e)$ the set of ascents of $e$ and by $\asc(e)$ the cardinality of $\Asc(e)$. An $n$-inversion sequence $e$ is called a \textit{derangement} if $e_n \neq 0$ and there exists no $i$ such that $e_i = e_{i+1} = 0$. We denote by $\D_n$ the set of derangements. 
\end{madef}
One has the following proposition.
\begin{prop}[\cite{Savage_2012} Lemma 1, \cite{Gustafsson_2020} Theorem 4.1]
For all $n \geq 2$ we have the two equalities 
\begin{align*}
    A_n(t) &= \sum_{e \in \I_n}t^{\asc(e)}, \\
    d_n(t) &= \sum_{e \in \D_n}t^{\asc(e)}.
\end{align*}
\end{prop}
The first equality comes from the bijection $\Le: \S_{[n]} \xrightarrow{\sim} \I_n$ (called ``Lehmer's code'') which sends a permutation $\sigma$ to the inversion sequence $(\Card(\{j \, | \,  j < i, \sigma(j) > \sigma(i)\}))_i$. This bijection also makes sense if we replace $[n]$ by any finite linearly ordered set of cardinality $n$. We have the following proposition. 
\begin{prop}[\cite{Hameister_2018} Theorem 5.1]
For all $n \geq 2$, one has the equality  
\begin{equation*}
    t\H(\U_{n-1, [n]})(t) = d_n(t).
\end{equation*}
\end{prop}
By Proposition \ref{propidealchowpoly} this is equivalent to the equality $\Hilb((h_{[n]})) = d_n(t)$, where  $h_{[n]} \in \CH(\B_{[n]})$. One of the main goals of this section is to reprove this result using Gröbner methods. To summarize, we will show that the normal monomials obtained after rewriting the elements $h_{[n]}m \,\,(m \in \NM([n]))$ are in bijection with $\D_n$ (via $\Phi$ and $\Le$). Unfortunately, the Gröbner basis we have used thus far is not well suited for this purpose. We will replace this Gröbner basis by another one with better rewriting, associated to a new presentation which we describe now. In this article we will only use this presentation for the Chow rings of boolean matroids but we believe this presentation may be of independent interest, and so we introduce it for arbitrary matroids. \\

We start by recalling the main result of \cite{Braden_2022} and its variation appearing in \cite{Pagaria_2023}. We refer to those articles for more details. For any matroid $\M$ on some ground set $E$ and $i$ some element of $E$, we define $\M\setminus i$ as the matroid on $E\setminus i$ with flats $F\setminus \{i\}$, where $F$ ranges over the flats of $\M$. For instance one has $\B_{E}\setminus i = \B_{E\setminus i}$. For any flat $F \in \L(\M)$ we define $\M^F$ as the matroid on the ground set $F$ with flats the flats of $\M$ contained in $F$, and $\M_F$ the matroid on the ground set $E\setminus F$ with flats the sets of the form $F'\setminus F$ where $F'$ is a flat of $\M$ containing $F$. For instance one has $(\B_{E})^F = \B_F$ and $(\B_{E})_F = \B_{E\setminus F}$. There is an injective morphism of rings $\theta_i: \CH(\M\setminus i) \rightarrow \CH(\M)$ defined by sending the generator $x_F$ to $x_F + x_{F \cup i}$ for all $F \in \L(\M) \setminus\{\emptyset\}$, with the convention that $x_S$ is zero if $S\subset E$ is not a flat of $\M$. We denote by $\CH_{(i)}(\M)$ the image of $\theta_i$ in $\CH(\M)$. Let $\Scol_i$ denote the set of flats $F$ of $M$ which contain $i$ and such that $F\setminus i$ is a nonempty flat of $\M$. 
\begin{theo}[\cite{Braden_2022} Theorem 1.1, \cite{Pagaria_2023} Theorem 5.4]\label{theosemismall}
There is a direct sum decomposition of $\CH(\M)$ into indecomposable graded $\CH(\M\setminus i)$-modules 
\begin{equation}\label{eqsemismalldec}
    \CH(\M) = \CH_{(i)}(\M)\oplus \bigoplus_{F \in \Scol_i}x_F\CH_{(i)}(\M).
\end{equation}
Moreover, for any $F \in \Scol_i$ we have an isomorphism 
\begin{equation*}
     \CH(\M^{F\setminus i})\otimes \CH(\M_F) \xrightarrow[\sim]{\psi^F} x_F\CH_{(i)}(\M).
\end{equation*}
where $\psi^{F}$ sends elements $\prod_{F ' \in V'}x_{F'} \otimes \prod_{F'' \in V''}x_{F''}$ to $x_F\prod_{F ' \in V'}x_{F'}\prod_{F'' \in V''}x_{F''\cup F}$, for any $V'$ and $V''$ collections of nonempty flats in $M^{F\setminus i}$ and $M_{F}$ respectively.
\end{theo}
This decomposition, commonly referred to as the semi-small decomposition, was first proved in \cite{Braden_2022}, with a small variation in the case where $E \in \Scol_i$. The form above comes from \cite{Pagaria_2023}. \\

Let $E$ be a finite linearly ordered set and $\M$ a matroid on $E$. For any subset $S = \{s_1 < \ldots < s_k \}\subset E$ we denote by $M\setminus S$ the matroid $\M\setminus i_k\setminus i_{k-1}\setminus ... \setminus i_1$, and by $\theta_{S}$ the morphism $\theta_{i_k}\circ \ldots \circ \theta_{i_1}$. We also denote $\theta = \theta_{\max E}$. 

\begin{madef}
For any nonempty flat $F$ we call the \textit{generation} of $F$ and denote by $\Gen(F)$ the smallest element of $E$ such that $F_{\leq \Gen(F)}$ generates $F$, that is, such that the intersection of the flats containing $F_{\leq \Gen(F)}$ is $F$. We put $$\s_F \coloneqq \theta_{E_{>\Gen(F)}}(x_{F_{\leq \Gen(F)}}),$$ where $x_{F_{\leq \Gen(F)}}$ is viewed as a generator in $\CH(\M\setminus E_{>\Gen(F)})$. 
\end{madef}
\begin{rmk}
A flat of $\M$ belongs to $\Scol_{\max E}$ if and only if its generation is $\max E$. More generally, a flat $F$ has generation $i$ if and only if $F_{\leq i}$ belongs to $\Scol_{i} \subset \L(\M\setminus E_{>i})$. 
\end{rmk}
\begin{rmk}\label{rmkinterpolation}
By using the explicit definition of $\theta$ we get 
\begin{equation*}
    \s_F = \sum_{S \subset E_{>\Gen(F)}} x_{F_{\leq \Gen(F)} \cup S} = \sum_{S \subset E_{>\Gen(F)}} x_{F \cup S}
\end{equation*}
with the second equality coming from the fact that since $x_G$ is zero when $G$ is not a flat, the only nonzero elements in the first sum are indexed by sets containing the whole $F$. If a flat $F$ has generation $\max E$ then we have $\s_F = x_F$. On the other hand, if $F$ is of the form $E_{\leq k}$ for some $k$, then we have $\s_F = h_F$. In other words, one could say that the generators $\s$ interpolate between the generators $x$ and the generators $h$ (via the linear order on $E$).  
\end{rmk}
Let us denote by $\L^{>1}(\M)$ the set of flats of $\M$ which are not the empty set and which are not of the form $\overline{\imath}$ for some $i \in E$. 
\begin{lemma}
    The elements $\s_F$, with $F$ running over $\L^{>1}(\M)$, generate $\CH(\M)$ as a ring.
\end{lemma}
\begin{proof}
We will prove this result by induction on the cardinality of $E$. For precision we put the matroid in which we are considering the generator in superscript. By construction of the generators $\s$ we have 
\begin{multline}\label{eqsgentransi}
    \{\s^\M_F, F \in \L^{>1}(\M)\} = \{x_F \, |\, F \in \Scol_{\max E}\} \, \cup
    \{\theta(\s^{\M\setminus \max E}_{F\setminus \max E}) \, | \, F \in \L^{>1}(\M), \Gen(F) < \max E\}.
\end{multline}
The semi-small decomposition \ref{theosemismall} states that every element $f \in \CH(\M)$ can be expressed as a sum
\begin{equation*}
    f = \theta(f') + \sum_{F \in \Scol_{\max E}}x_{F}\theta(f_F)
\end{equation*}
for some $f', f_F \, (F \in \Scol_{\max E}) \in \CH(\M\setminus \max E)$. Equality \eqref{eqsgentransi} and our induction hypothesis conclude the proof. 
\end{proof}
We denote by $\Is(\M)$ the kernel of the natural morphism $\Q[\s_F, F\in \L^{>1}(\M)] \rightarrow \CH(\M)$. We turn our attention back to boolean matroids. The next proposition gives a Gröbner basis of $\Is(\B_{E})$ for any finite linearly ordered set $E$. 
\begin{prop}\label{propgbasiss}
Let $(E, <)$ be any finite linearly ordered set. Let $\vartriangleleft$ be any total order on $\L^{>1}(\B_{E})$ such that for any two flats $F_1, F_2 \in  \L^{>1}(\B_{E})$, if $\Gen(F_1) < \Gen(F_2)$ then $F_1 \vartriangleleft F_2$, and if $\Gen(F_1) = \Gen(F_2)$ and $F_1 \supset F_2$ then $F_1 \vartriangleleft F_2$. The degree $2$ part $\Isq(\B_{E})$ of $\Is(\B_{E})$ forms a quadratic Gröbner basis of $\Is(\B_{E})$ for the degree lexicographic ordering on the monomials induced by $\vartriangleleft$. In addition, a quadratic monomial $\s_{S_1}\s_{S_2}$ is normal with respect to that Gröbner basis if and only if $h_{S_1}h_{S_2}$ is a normal monomial of the Gröbner basis described in Proposition \ref{propgbasish}.
\end{prop}
\begin{proof}
By Proposition \ref{propgbasish} and Proposition \ref{propmaingrob} it is enough to prove the second part of the statement, and for this second part it is enough to prove that if $\s_{S_1}\s_{S_2}$ is normal with respect to $\Isq(\B_E)$ then $h_{S_1}h_{S_2}$ belongs to the set $\eqref{eqnormal}.$ Assume $\Gen(S_1) \geq \Gen(S_2).$ Note that for boolean matroids, the generation of a flat is simply its maximum. By construction of the generators $\s$ we have $\s_{S_1} = \theta_{E_{>\Gen(S_1)}}(\s^{\M\setminus E_{> \Gen(S_1)}}_{S_1})$ and $\s_{S_2} = \theta_{E_{>\Gen(S_1)}}(\s^{\M\setminus E_{> \Gen(S_1)}}_{S_2})$, where the matroid in superscript indicates in which matroid we are considering the generator. Since $\theta$ is injective and preserves the order on the generators we can assume $\Gen(S_1) = \max E$. If $S_1$ contains $S_2$ it is enough to prove $\Gen(S_1) > \Gen(S_2).$ If we had $\Gen(S_1) = \Gen(S_2) (=\max E)$, then by decomposition \eqref{eqsemismalldec} we would have an equation of the form 
\begin{equation}\label{eqsamegeneration}
    \s_{S_1}\s_{S_2} = \theta(f) + \sum_{F\in \Scol_{\max E}}\s_F\theta(f_F)
\end{equation}
for some elements $f, f_F \, (F \in \Scol_{\max E})$ in $\CH(\M\setminus \max E)$. By definition of the generators $\s$ we can express $\theta(f)$ and the $\theta(f_F)$'s using generators $\s$ of smaller generations. By the definition of our order on the generators $\s$ this implies that $\s_{S_1}\s_{S_2}$ is the leading term of equation \eqref{eqsamegeneration}, and so $\s_{S_1}\s_{S_2}$ cannot be normal. If $S_2$ is not contained in $S_1$ we have 
\begin{align*}
    \s_{S_1}\s_{S_2} &=  x_{S_1}\biggl(\sum_{S \subset E_{>\Gen(S_2)}}x_{S_2 \cup S}\biggr) \\
                   &= \sum_{S \subset E_{>\Gen(S_2)}}x_{S_1}x_{S_2 \cup S}
\end{align*}
and since the product of two $x$-generators is zero unless one flat is included in the other, we see that the above sum is zero unless $S_2$ contains every element of $S_1$ below $\Gen(S_2)$, that is, unless $S_2$ is an initial segment of $S_1 \cup S_2$. Denote by $i$ the biggest element of $S_2$ which is not in $S_1$. Denote also $I \coloneqq (S_2)_{\leq i}$. This is the smallest initial segment of $S_1 \cup S_2$ such that $S_1 \cup I = S_1 \cup S_2$. One has 
\begin{align*}
    \s_{S_1}\s_{I} &= \sum_{S \subset E_{>i}}x_{S_1}x_{I \cup S} \\
                 &= \sum_{\substack{S \subset E_{>i} \\ S \supset (S_1)_{\geq i}}}x_{S_1}x_{I \cup S} \\
                 &= \s_{S_1}\s_{S_2}.
\end{align*}
By definition of our order on the generators $\s$ we have $I \trianglelefteq S_2$ and so if $\s_{S_1}\s_{S_2}$ is normal we must have $S_2 = I$. What is left to prove is that we have $\Card(S_2 \setminus S_1) > 1$, or equivalently that $S_2$ cannot be equal to  $(S_1)_{< i}\cup \{i\}$. We have 
\begin{align*}
    h_{i} & = 0 \\
    & = \sum_{\substack{S \subset E \\ S \ni i}} x_{S} \\
    & = \sum_{S' \subset E_{<i}}\sum_{S'' \subset E'_{>i}}x_{S'\cup S'' \cup \{i\}} \\
    & = \sum_{S' \subset E_{<i}}\s_{S'\cup \{i\}}.
\end{align*}
Applying $\s_{S_1}$ we get 
\begin{align*}
    0 & = \s_{S_1}\biggl(\sum_{S' \subset E_{<i}}\s_{S'\cup \{i\}}\biggr) \\
    & = \sum_{\substack{S' \subset E_{<i} \\ S' \supset (S_1)_{< i}} }\s_{S_1}\s_{S'\cup \{i\}},
\end{align*}
and the leading term of that relation is $\s_{S_1}\s_{(S_1)_{< i}\cup \{i\}}$.
\end{proof}
We denote by $\NM_{\ss}(E)$ the normal monomials associated to the Gröbner basis described above. There is an obvious bijection $\NM(E) \xrightarrow[\sim]{f} \NM_{\ss}(E)$ given by replacing $h$ by $\s$. We have the following lemma which describes how $\psi^F$ interacts with the generators $\s$. 
\begin{lemma}\label{lemmapsi}
    For any $F \in \Scol_{\max E}$ different from $E$, and any monomial $m = \prod_{S \in V}\s_{S}$ in $\CH(\B_{E\setminus F})$ with $V$ some set of subsets of $E\setminus F$, we have
    \begin{equation*}
        \psi^F(1 \otimes m) = \s_F \prod_{S \in V} \s_{S \cup F_{\leq \max S}}.
    \end{equation*}
\end{lemma}
\begin{proof}
        One has 
        \begin{align*}
            \psi^F(1\otimes m) &= \psi^F\left(1 \otimes \prod_{S\in V}\left(\sum_{\substack{S' \subset (E\setminus F)_{> \max S}}}x_{S\cup S'}\right)\right) \\
            &= \s_F \prod_{S \in V}\left(\sum_{\substack{S' \subset (E\setminus F)_{> \max S}}} x_{S\cup S' \cup F} \right)\\
            &= \s_F \prod_{S \in V} \left(\sum_{\substack{S' \subset E_{> \max S}}} x_{S\cup F_{\leq \max S} \cup S'}\right) \\
            &= \s_F \prod_{S \in V}\s_{S \cup F_{\leq \max S}}.
        \end{align*}
\end{proof}
We now describe the rewriting of terms of the form $\s_{E}m$ where $m$ is a normal monomial in $\NM_{\ss}(E)$. A priori this rewriting could be a sum of a several normal monomials, but crucially, in our case the rewriting will yield zero or a single normal monomial. Our philosophy is that this rewriting should have a geometric interpretation in terms of inversion sequences. This means we should try to describe the map $\Le(\Phi(f^{-1}(\s_{E}f(\Psi(\Le^{-1}(-))))))$. In the name of readability we will simply denote this map by ``$\s_{n}(-)$'' with $n$ the cardinality of $E$.
\begin{prop}\label{propcomputation}
Let $n\geq 2$ be an integer, and $E$ a finite linearly ordered set of cardinality $n$. For any normal monomial $m \in \NM_{\ss}(E)$, the term $\s_E m$ rewrites to a normal monomial or to zero. For any inversion sequence $e \in \I_n$, the inversion sequence $\s_n(e)$ can be computed inductively as follow : 
\begin{enumerate}
\item If $e$ has no zero after $1$ then $\s_n (e)$ is the inversion sequence obtained by applying $\s_{n-1}$ to the inversion sequence $(e_i)_{i \geq 2}$ viewed as an $n-1$-inversion sequence. More precisely, if we set $\Tilde{e}_i \coloneqq e_{i+1} - 1$ for $i \leq n - 1$ then we have $\s_n(e)_1 = 0$ and $\s_n(e)_k = (\s_{n-1}(\Tilde{e}))_{k-1} + 1$ for all $2\leq k \leq n$.
\item If $e_n = 0$ then $\s_n (e)$ is the inversion sequence obtained from $e$ by removing the $n$-th entry, moving the graph of $e$ of one unit in direction $(1, 1)$ and setting the first entry to $0$. More precisely, we have $(\s_n (e))_1 = 0$ and $\s_n(e)_k = e_{k-1} + 1$ for all $2 \leq k \leq n$.
\item Finally, if $e$ has a zero after $1$ and its maximal zero is strictly less than $n$, let us denote by $k$ this maximal zero. The inversion sequence $\s_n(e)$ is the inversion sequence obtained from $e$ by keeping the same values after $k$ and applying $\s_{k-1}$ to the inversion sequence before $k$. More precisely, if we denote by $\Tilde{e}$ the $k-1$ inversion sequence $(e_i)_{i\leq k-1}$ we have $\s_n(e)_j = \s_{k-1}(\Tilde{e})_j$ for all $1 \leq j \leq k-1$ and we have $\s_n(e)_j = e_j$ for all $k \leq j \leq n$.  
\end{enumerate}
\end{prop}
Here are some examples of rewritings of inversion sequences below. 
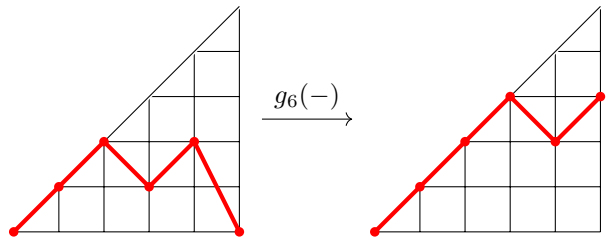
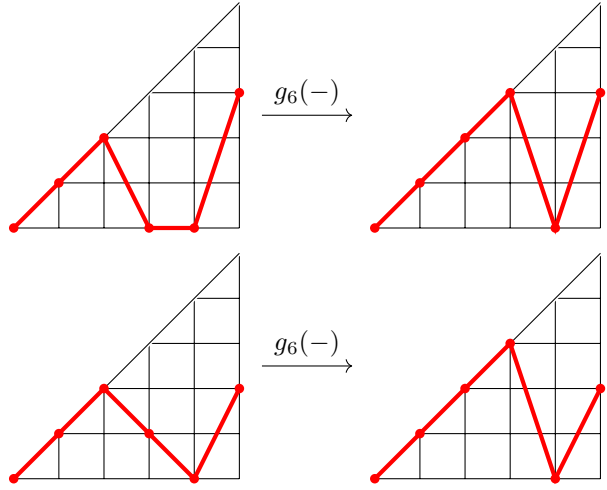
\begin{ex}
\begin{equation*}
\begin{tikzpicture}[scale=0.60, baseline=14px]
    \draw[step=1.0,black,thin] (0,0) grid (5,5);
    \draw [ultra thick, draw=white, fill=white]
       (0,0) -- (5,5) -- (0,5) -- cycle;
    \draw[black] (0,0) -- (5,5);
    \draw[ultra thick, red] (0,0) -- (2,2) -- (3,1) -- (4,2) -- (5,0);
    \fill[red] (0,0) circle (3pt);
    \fill[red] (1,1) circle (3pt);
    \fill[red] (2,2) circle (3pt);
    \fill[red] (3,1) circle (3pt);
    \fill[red] (4,2) circle (3pt);
    \fill[red] (5,0) circle (3pt);
    \draw[->] (5.5, 2.5) -- (7.5, 2.5) node[midway,above] {$\s_6(-)$};

    \draw[step=1.0,black,thin] (8,0) grid (13,5);
    \draw [ultra thick, draw=white, fill=white]
       (8,0) -- (13,5) -- (8,5) -- cycle;
    \draw[black] (8,0) -- (13,5);
    \draw[ultra thick, red] (8,0) -- (9,1) -- (10,2) -- (11,3) -- (12,2) -- (13,3);
    \fill[red] (8,0) circle (3pt);
    \fill[red] (9,1) circle (3pt);
    \fill[red] (10,2) circle (3pt);
    \fill[red] (11,3) circle (3pt);
    \fill[red] (12,2) circle (3pt);
    \fill[red] (13,3) circle (3pt);
    
\end{tikzpicture}
\end{equation*}
\begin{equation*}
\begin{tikzpicture}[scale=0.60, baseline=14px]
    \draw[step=1.0,black,thin] (0,0) grid (5,5);
    \draw [ultra thick, draw=white, fill=white]
       (0,0) -- (5,5) -- (0,5) -- cycle;
    \draw[black] (0,0) -- (5,5);
    \draw[ultra thick, red] (0,0) -- (2,2) -- (3,0) -- (4,0) -- (5,3);
    \fill[red] (0,0) circle (3pt);
    \fill[red] (1,1) circle (3pt);
    \fill[red] (2,2) circle (3pt);
    \fill[red] (3,0) circle (3pt);
    \fill[red] (4,0) circle (3pt);
    \fill[red] (5,3) circle (3pt);
    \draw[->] (5.5, 2.5) -- (7.5, 2.5) node[midway,above] {$\s_6(-)$};

    \draw[step=1.0,black,thin] (8,0) grid (13,5);
    \draw [ultra thick, draw=white, fill=white]
       (8,0) -- (13,5) -- (8,5) -- cycle;
    \draw[black] (8,0) -- (13,5);
    \draw[ultra thick, red] (8,0) -- (9,1) -- (10,2) -- (11,3) -- (12,0) -- (13,3);
    \fill[red] (8,0) circle (3pt);
    \fill[red] (9,1) circle (3pt);
    \fill[red] (10,2) circle (3pt);
    \fill[red] (11,3) circle (3pt);
    \fill[red] (12,0) circle (3pt);
    \fill[red] (13,3) circle (3pt);
    
\end{tikzpicture}
\end{equation*}  
\begin{equation*}
\begin{tikzpicture}[scale=0.60, baseline=14px]
    \draw[step=1.0,black,thin] (0,0) grid (5,5);
    \draw [ultra thick, draw=white, fill=white]
       (0,0) -- (5,5) -- (0,5) -- cycle;
    \draw[black] (0,0) -- (5,5);
    \draw[ultra thick, red] (0,0) -- (2,2) -- (3,1) -- (4,0) -- (5,2);
    \fill[red] (0,0) circle (3pt);
    \fill[red] (1,1) circle (3pt);
    \fill[red] (2,2) circle (3pt);
    \fill[red] (3,1) circle (3pt);
    \fill[red] (4,0) circle (3pt);
    \fill[red] (5,2) circle (3pt);
    \draw[->] (5.5, 2.5) -- (7.5, 2.5) node[midway,above] {$\s_6(-)$};

    \draw[step=1.0,black,thin] (8,0) grid (13,5);
    \draw [ultra thick, draw=white, fill=white]
       (8,0) -- (13,5) -- (8,5) -- cycle;
    \draw[black] (8,0) -- (13,5);
    \draw[ultra thick, red] (8,0) -- (9,1) -- (10,2) -- (11,3) -- (12,0) -- (13,2);
    \fill[red] (8,0) circle (3pt);
    \fill[red] (9,1) circle (3pt);
    \fill[red] (10,2) circle (3pt);
    \fill[red] (11,3) circle (3pt);
    \fill[red] (12,0) circle (3pt);
    \fill[red] (13,2) circle (3pt);
    
\end{tikzpicture}
\end{equation*}  
\end{ex}
\begin{proof}
Let $n \geq 2$ be an integer and $E$ some finite linearly ordered set of cardinality $n$. The first statement will come by the proof of computations 1), 2), 3). 
\begin{enumerate}
     \item If the inversion sequence $e$ satisfies $e_n = 0$ then this means that we have $\Le^{-1}(e)(n) = \max E.$ This means that the monomial $f(\Psi(\Le^{-1}(e)))$ is of the form $\theta(m')$ with $m'$ a normal monomial in $\CH(\B_{E\setminus \max E})$, whose associated inversion sequence is simply $(e_i)_{i\leq n-1}$. By the description of the quadratic normal monomials \eqref{eqnormal}, applying $\s_E$ to a normal monomial $\theta(m')$ preserves normality, and, by the formula defining $\Psi$, the inversion sequence associated to $\s_E\theta(m')$ is  obtained from that of $m'$ by shifting it in direction $(1,1)$, as described in computation 1). 
    \item If the inversion sequence $e$ has no zero after $1$ then this means that we have $\Le^{-1}(e)(1) = n.$ This means that the monomial $f(\Psi(\Le^{-1}(e)))$ is of the form $\s_{E}\theta(m')$ with $m'$ a normal monomial in $\CH(\B_{E\setminus \max E})$ (whose associated inversion sequence is $\Tilde{e}$ described in computation 2). We then have 
    \begin{align*}
    \s_E f(\Psi(\Le^{-1}(e))) &= \s_E \s_{E}\theta(m') \\
                                 &= \s_E \s_{E\setminus \max E} \theta(m') \\
                                 &= \s_E \theta(\s^{E\setminus \max E}_{E\setminus \max E}m').
    \end{align*} 
    The second equality comes from the identity $\s_E^2 = (h_E)^2 = h_{E}h_{E\setminus \max E} = \s_E \s_{E\setminus \max E}$. The inversion sequence associated to the rewriting of $\s^{E\setminus \max E}_{E\setminus \max E}m'$ is $\s_{n-1}\Tilde{e}$. Furthermore, as noted previously the injective morphism $\theta$ preserves the order on the $\s$-generators and so it commutes with the rewriting. Finally, as observed in the proof of computation 1), applying $\s_E$ has the desired effect on the inversion sequence of $\s_{n-1}\Tilde{e}$. 
    \item If $e$ has a zero after $1$ and its maximal zero is not $n$, denote by $k$ the maximal zero of $e$ and denote by $\Tilde{e}$ the inversion sequence $(e_i)_{i \leq k-1}$. Denote also $\sigma \coloneqq \Le^{-1}(e)$. Since $k$ is the maximal zero of $e$ we must have $\sigma(k) = n$, and $k$ is a descent of $\sigma$. Let us denote $F = S(\sigma, k).$ The normal monomial associated to $e$ is of the form $\s_F m$ with $m$ a normal monomial in $\theta(\CH(\B_{E \setminus \max E}))$. Let us denote by $V_{\bot F}$ the set of subsets $S$ of $E\setminus \max E$ indexing a generator of $m$ and such that $S$ is not comparable with $F$. Denote also by $V_{< F}$ the set of subsets $S$ of $E\setminus \max E $ indexing a generator in $m$ and such that $S$ is contained in $F$, so that we have $m = \prod_{S \in V_{\bot F}}\s_S\prod_{S' \in V_{< F}} \s_{S'}$. Since $\s_F m$ is normal, for every $S$ in $V_{\bot F}$ the subset $S$ must be an initial segment of $S\cup F$. This implies that we have $(\psi^F)^{-1}(\s_F \prod_{S \in V_{\bot F}}\s_S ) = 1 \otimes \prod_{S \in V_{\bot F}}\s_{S\setminus F}$.
    \begin{lemma}\label{lemmabeginning}
    The monomial $\Psi(\Le^{-1}(\Tilde{e})) \in \NM_{\ss}(E \setminus F)$ is $\prod_{S \in V_{\bot F}}\s_{S\setminus F}$.
    \end{lemma}
    \begin{proof}
    The bijection $\Le^{-1}(\Tilde{e})$ is the restriction $\sigma_{|[k-1]}$. For all descents $i$ of $\sigma_{|[k-1]}$ one has $S(\sigma_{|[k-1]}, i) = S(\sigma, i)\setminus F$ and the sets $S(\sigma, i), i \leq k-1$ are exactly the elements of $V_{\bot G}$.
    \end{proof}
    In the proof of Proposition \ref{propgbasiss} we have proved the identity 
    \begin{equation*}
        \s_E \s_F = \s_F \s_{E_{\leq \max E\setminus F}}.
    \end{equation*}
    This gives
    \begin{align*}
        \s_E\s_F m & = \s_F \s_{E_{\leq \max E\setminus F}}\prod_{S \in V_{\bot F}}\s_S\prod_{S' \in V_{< F}} \s_{S'} \\
        &= \psi^F((\psi^F)^{-1}( \s_F \s_{E_{\leq \max E\setminus F}}\prod_{S \in V_{\bot F}}\s_{S}))\prod_{S \in V_{< F}} \s_{S} \\
        &= \psi^F(1\otimes \s_{E\setminus F}\prod_{S \in V_{\bot F}}\s_{S\setminus F})\prod_{S \in V_{< F}} \s_{S}.
    \end{align*}
    Denote by $m' = \prod_{S' \in V'}\s_{S'}$ the rewriting of the monomial $\s_{E\setminus F}\prod_{S \in V_{\bot F}}\s_{S\setminus F}$ in $\CH(\B_{E\setminus F})$, provided this rewriting is not zero. Denote by $\sigma'$ the bijection associated to this normal monomial.  By Lemma \ref{lemmapsi} we get
    \begin{equation}\label{eqlastmonomial}
        \s_E\s_F m = \s_F \prod_{S' \in V'}\s_{S' \cup F_{\leq \max S}}\prod_{S \in V_{< F}} \s_{S}.
    \end{equation}
    Denote $\sigma'' \coloneqq \sigma'(1)\ldots \sigma'(k-1) \sigma(k) \ldots \sigma(n)$ and $e'' = \Le(\sigma'').$
    \begin{lemma}
        The monomial associated to $\sigma''$ is the right-hand side of \eqref{eqlastmonomial}, and the inversion sequence $e''$ is the inversion sequence described in computation 3).
    \end{lemma}
    \begin{proof}
    For every descent $i \geq k$ of $\sigma''$ one has $S(\sigma'', i) = S(\sigma, i)$ and those sets are exactly $\{F\}\cup \{S \in V_{< F}\}.$ For every descent $i \leq k-1$ of $\sigma''$ one has 
    \begin{align*}
         S(\sigma'', i) &= S(\sigma', i) \cup \{j \, | \, j \geq k, \sigma(j) \leq i\} \\
         &= S(\sigma', i) \cup F_{\leq i} \\
         &= S(\sigma', i) \cup F_{\leq \max  S(\sigma', i)},
    \end{align*}
    and those sets are exactly $\{S \cup F_{\leq \max S}, S \in V'\}.$
    For every $i \geq k$ one has 
     \begin{align*}
        e'_{i} &= \Card(\{j \, | \, j < i, j\geq k, \sigma(j) > \sigma(i)\}) + \Card(\{j \, | \, j < k, \sigma'(j) > \sigma(i)\}) \\
        &= \Card(\{j \, | \, j < i, j\geq k, \sigma(j) > \sigma(i)\}) + \Card((E\setminus F)_{\geq \sigma(i)}) \\
        &= \Card(\{j \, | \, j < i, j\geq k, \sigma(j) > \sigma(i)\}) + \Card(\{j \, | \, j < k, \sigma(j) > \sigma(i)\}) \\
        &= e_{i},
    \end{align*} 
    and by Lemma \ref{lemmabeginning} we have $(e''_i)_{i\leq k-1} = \s_{k-1}\Tilde{e}$.
    \end{proof}
    By Proposition \ref{propbij} the right-hand side of \eqref{eqlastmonomial} must be the rewriting of $\s_E\s_F m$ which concludes the proof.
\end{enumerate}    
\end{proof}
For all $n \geq 2$ and $e \in \I_n$ an inversion sequence such that $e_2 =1$, we denote by $\fp(e)$ (standing for ``first peak'') the inversion sequence defined by 
$\fp(e)_i = e_{i+1} - 1$ for all $1\leq i \leq \fz(e)-2$ where $\fz(e)$ is the first zero of $e_i$ after $1$ if such a zero exists, and $n+1$ otherwise. For all $n \geq 1$ and $1 \leq k \leq n$ we define a set $\D^k_n$ of $n$-inversion sequences as follow. For all $n \geq 1$ we put $\D^1_n \coloneqq \D_n$, with the convention $\D^1_1 = \emptyset$. Then, for all $k\leq n$ we define $\D^k_n$ as the set of inversion sequences $e \in \I_n$ such that $e$ belongs to $\D^1_n$ and $\fp(e)$ belongs to $\D^{k-1}_{\fz(e) - 2}$. 
\begin{coro}
For all $n\geq 2$ and $1\leq k\leq n-1$ we have
\begin{equation*}
    t^k\H(\U_{n-k, n}) = \sum_{e \in \D^k_n}t^{\asc(e)}.
\end{equation*}
\end{coro} 
\begin{proof}
For all $n\geq 2$ and $1\leq k \leq n-1$, let us denote by $(\D^k_{n})'$ the set inversion sequences associated to normal monomials which can be obtained as the rewriting of terms of the form $\s_n^k m$, with $m$ a normal monomial. By Proposition \ref{propcomputation} those monomials span the ideal $(\s_n^k)$, and since normal monomials are linearly independent they form a basis of that ideal. By Proposition \ref{propidealchowpoly} we get 
\begin{equation*}
    t^k\H(\U_{n-k, n}) = \sum_{e \in (\D^k_n)'}t^{\asc(e)}.
\end{equation*}
We need to prove the equality $(\D^k_n)' = \D^k_n$. We first prove the inclusion $(\D^k_n)' \subset \D^k_n$ by induction. By Proposition \ref{propcomputation} it is clear that we have $(\D^1_n)' \subset \D^1_n$. Let $m$ be the rewriting of some term of the form $\s_n^k m'$ for some normal monomial $m'$ and $k\geq 2$. Denote by $m''$ the rewriting of $\s^{k-1}_n m'$, and denote by $e''$ its associated inversion sequence. By our induction hypothesis we have $e'' \in \D^{k-1}_n$. The monomial $m$ is the rewriting of $\s_n m''$ and so by Proposition \ref{propcomputation} (more precisely iterated use of computation 3) and then computation 1)) its associated inversion sequence is given by $1 + \s_{\fz(e'') -2}\fp(e'')$ before $\fz(e'')$ and by $e''$ after $\fz(e'')$. By our induction hypothesis this means that the inversion sequence associated to $m$ belongs to $\D^k_n$. We now prove the inclusion $\D^k_n \subset (\D^k_n)'$, also by induction. For $k=1$, let $e$ be an inversion sequence in $\D^1_n.$ By definition of $\fz(e)$ the inversion sequence $(e_i)_{i < \fz(e)}$ has no zero after $1$, and so it is the shift in direction $(1,1)$ of some inversion sequence $(e'_i)_{i < \fz(e) - 1}$. By Proposition \ref{propcomputation} we have that $e''$ is $\s_n$ applied to the inversion sequence given by $(e'_i)$ before $\fz(e) - 1$, $0$ at $\fz(e)-1$, and $e$ after $\fz(e)-1$. For $k \geq 2$, let $e$ be an inversion sequence in $\D^k_n$. By our induction hypothesis $\fp(e)$ belongs to $(\D^{k-1}_{\fz(e)-2})'$, that is, is of the form $\s_{\fz(e)-2}(e')$ for some $e' \in \D^{k-2}_{\fz(e)-2}$. By Proposition \ref{propcomputation}, the inversion sequence $e$ can be obtained as the rewriting $\s_n(e'')$ with $e''$ the inversion sequence given by $e' + 1$ before $\fz(e) -1$ and $e$ after. The inversion sequence $e''$ belongs to $\D^{k-1}_n$, which is equal to $(\D^{k-1}_n)'$ by our induction hypothesis. This means that we have $\s_n(e'') \in (\D^{k}_n)'$ which concludes the proof. 
\end{proof}
\begin{rmk}
    Unfortunately, the proofs of the real-rootedness of $\H(\B_{[n]})$ and $\H(\U_{n-1, n})$ described in the introduction do not generalize verbatim in higher corank because if we put $d^{k,i}_n(t) \coloneqq \sum\limits_{\substack{e\in \D^k_n \\ e_n = i}} t^{\asc(e)}$ for $i \geq 1$ and $\d_n^{k,0} \coloneqq \sum_{e \in \D^k_{n-1}}$, computations by computer show that the sequence of polynomials $(d^{k,i}_n)_i$ is not interlacing. However, computations also show that some small variations around that sequence of polynomials are interlacing. For instance, for every $n \leq 25$ the sequence $(d^{2,0}_n, d^{2, 1}_n, d^{2,3}_n, \ldots)$ is interlacing, and so is the sequence $(d^{2,0}_n, d^{2,1}_n + d^{2,2}_n, d^{2,3}_n, \ldots)$. This seems to indicate that one would ``just'' need to find the right variation in order to set up the induction and prove the real-rootedness of $\H(\U_{n-k, n})$ for $k\geq 2.$
\end{rmk}
    
\bibliographystyle{alpha}
\bibliography{sample}

\newcommand{\etalchar}[1]{$^{#1}$}
\begin{thebibliography}{BHM{\etalchar{+}}22}

\bibitem[AHK18]{Adiprasito_2018}
Kareem Adiprasito, June Huh, and Eric Katz.
\newblock {Hodge theory for combinatorial geometries}.
\newblock {\em Annals of Mathematics}, 188(2):381 -- 452, 2018.

\bibitem[BES20]{Backman_Spencer_Eur_2020}
Spencer Backman, Christopher Eur, and Connor Simpson.
\newblock Simplicial generation of {C}how rings of matroids.
\newblock {\em S\'{e}m. Lothar. Combin.}, 84B:Art. 52, 11, 2020.

\bibitem[BHM{\etalchar{+}}22]{Braden_2022}
Tom Braden, June Huh, Jacob~P. Matherne, Nicholas Proudfoot, and Botong Wang.
\newblock A semi-small decomposition of the {C}how ring of a matroid.
\newblock {\em Advances in Mathematics}, 409:108646, 2022.

\bibitem[Bon15]{Bona_2015}
Miklos Bona.
\newblock {\em Handbook of enumerative combinatorics}.
\newblock Chapman and Hall, 2015.

\bibitem[BW93]{BW_1993}
Thomas Becker and Volker Weispfenning.
\newblock {\em Gröbner Bases, A Computational Approach to Commutative Algebra}.
\newblock Springer, 1993.

\bibitem[Cor23]{coron_2023}
Basile Coron.
\newblock Supersolvability of built lattices and {K}oszulness of generalized {C}how rings, 2023.

\bibitem[Dan78]{Danilov_1978}
V~I Danilov.
\newblock The geometry of toric varieties.
\newblock {\em Russian Mathematical Surveys}, 33(2):97, apr 1978.

\bibitem[DCP95]{de_concini_wonderful_1995}
Corrado De~Concini and Claudio Procesi.
\newblock Wonderful models of subspace arrangements.
\newblock {\em Selecta Mathematica. New Series}, 1(3):459--494, 1995.

\bibitem[FMSV24]{FERRONI_2024}
Luis Ferroni, Jacob~P. Matherne, Matthew Stevens, and Lorenzo Vecchi.
\newblock Hilbert–{P}oincaré series of matroid {C}how rings and intersection cohomology.
\newblock {\em Advances in Mathematics}, 449:109733, 2024.

\bibitem[Fro10]{Frobenius_1910}
G.~Frobenius.
\newblock {\em {\"U}ber die Bernoullischen Zahlen und die Eulerschen Polynome}.
\newblock Preussische Akademie der Wissenschaften Berlin: Sitzungsberichte der Preussischen Akademie der Wissenschaften zu Berlin. Reichsdr., 1910.

\bibitem[FY04]{feichtner_chow_2003}
Eva~Maria Feichtner and Sergey Yuzvinsky.
\newblock Chow rings of toric varieties defined by atomic lattices.
\newblock {\em Invent. Math.}, 155(3):515--536, 2004.

\bibitem[GS20]{Gustafsson_2020}
Nils Gustafsson and Liam Solus.
\newblock Derangements, {E}hrhart theory, and local h-polynomials.
\newblock {\em Advances in Mathematics}, 369:107169, 2020.

\bibitem[HRS18]{Hameister_2018}
Thomas Hameister, Sujit Rao, and Connor Simpson.
\newblock Chow rings of vector space matroids.
\newblock {\em Journal of Combinatorics}, 2018.

\bibitem[MM21]{Mastroeni_2021}
Matthew Mastroeni and Jason McCullough.
\newblock Chow rings of matroids are {K}oszul.
\newblock {\em Matematische Annalen}, 10 2021.

\bibitem[Pet15]{Pet_2015}
T.Kyle Petersen.
\newblock {\em Eulerian Numbers}.
\newblock Springer, 2015.

\bibitem[PP05]{Polischuk_2005}
Alexander Polishchuk and Leonid Positselski.
\newblock {\em Quadratic algebras}.
\newblock American Mathematical Society, 2005.

\bibitem[PP23]{Pagaria_2023}
Roberto Pagaria and Gian~Marco Pezzoli.
\newblock {Hodge Theory for Polymatroids}.
\newblock {\em International Mathematics Research Notices}, 2023(23):20118--20168, 03 2023.

\bibitem[SS12]{Savage_2012}
Carla~D. Savage and Michael~J. Schuster.
\newblock Ehrhart series of lecture hall polytopes and {E}ulerian polynomials for inversion sequences.
\newblock {\em Journal of Combinatorial Theory, Series A}, 119(4):850--870, 2012.

\bibitem[Sta72]{Stanley_1972}
R.~P. Stanley.
\newblock Supersolvable lattices.
\newblock {\em Algebra Universalis}, 2:197--217, 1972.

\bibitem[Wel76]{welsh_matroid_1976}
Dominic J.~A. Welsh.
\newblock {\em Matroid theory}.
\newblock L. {M}. {S}. {Monographs}, {No}. 8. Academic Press [Harcourt Brace Jovanovich, Publishers], London-New York, 1976.

\bibitem[Yuz02]{Yuzvinsky_2002}
Sergey Yuzvinsky.
\newblock Small rational model of subspace complement.
\newblock {\em Transactions of the American Mathematical Society}, 354(5):1921--1945, 2002.

\end{thebibliography}

\end{document}